\documentclass[12pt,reqno]{amsart}
\newtheorem{theorem}{Theorem}[section]

\newtheorem{corollary}[theorem]{Corollary}
\theoremstyle{definition}   

\newtheorem{example}[theorem]{Example}

\theoremstyle{remark}
\newtheorem{remark}[theorem]{Remark}

\numberwithin{equation}{section}

\usepackage{times}
\usepackage{enumerate}
\usepackage{graphicx,adjustbox}
\usepackage{hyperref}
\usepackage{amsmath,amssymb}
\usepackage{leftidx}

\usepackage{amscd}
\usepackage{graphicx}
\usepackage[all]{xy}
\title[Tangent Cones of Concatenated Numerical Semigroups  ]
{Tangent Cones of Bresinsky and Arsalan Curves}
\author{
Ranjana Mehta
\and
Joydip Saha
}
\date{}

\address{\small \rm  Department of Mathematics, SRM University AP, Amaravati 522240, Andhra Pradesh, India}
\email{ranjana.m@srmap.edu.in}

\address{\small \rm  Department of Mathematics, Barasat College, 1 Kalyani, Road, Barasat, West Bengal, Pin Code-700126, India.} 
\email{saha.joydip56@gmail.com}

\date{}



\allowdisplaybreaks

\begin{document}

\begin{abstract}
In this paper, we study the Ap\'{e}ry tables for the numerical semigroups given by Bresinsky and Arslan. Using the Ap\'{e}ry tables we write the tangent cones of the Bresinsky and Arsalan curves at the origin. Further, we calculate Hilbert series of the tangent cone of the Bresinsky and Arslan curves. We prove that both classes of the curve have Cohen-Macaulay tangent cone. 
\end{abstract}
\maketitle

\section{introduction}
Ap\'{e}ry table of a numerical semigroup associated to an affine monomial curve plays an important role in the characterizing invariants of its tangent cone. Using the Ap\'{e}ry table, we can calculate the explicit tangent cone and it's Hilbert series. We can also study the Cohen-Macauley and Buchbaum properties of the tangent cone; see the papers \cite{bjz}, \cite{cz1},\cite{cz2}.
\medskip

In \cite{ars}, F.Arslan studied the Cohen Macaulayness of the tangent cone of the Arslan curves using Gr\"{o}bner basis. F. Arslan further showed that in every affine $l-$ space with $l \geq  4$, there are monomial curves having a Cohen-Macaulay tangent cone with an arbitrarily large minimal number of generators which contradicts the case $l=3$, studied by Robbiano and Valla in \cite{LG}.
F. Arslan also determined the Hilbert series of the tangent cone of the Arslan curves and their extended versions. In \cite{hd} using Gr\"{o}bner basis, Herzog and Stamate showed that the tangent cone of the Bresinsky curves is Cohen–Macaulay.
\medskip

 In this paper, we compute the Ap\'{e}ry set and Ap\'{e}ry table for the numerical semigroups given by Bresinsky and Arslan. It came to our surprise that both classes have similar types of elements in the Ap\'{e}ry set and each element of the Ap\'{e}ry set has a unique expression. Using the Ap\'{e}ry tables, we write the tangent cones at the origin explicitly and further we calculate the Hilbert series of tangent cones. We prove that both classes have Cohen-Macaulay tangent cones.

\section{Preliminaries and Notations}
Let $\Gamma$ be a numerical semigroup. It is true that 
(see \cite{rgs}) the set $\mathbb{N}\setminus \Gamma$ is 
finite and that the semigroup $\Gamma$ has a unique minimal system of generators 
$ a_{1} < \cdots < a_{e}$. The greatest integer not belonging to $\Gamma$ 
is called the \textit{Frobenius number} of $\Gamma$, denoted by $F(\Gamma)$. The integers $a_{1}$ and $e$ are known as the \textit{multiplicity} and the 
\textit{embedding dimension} of the semigroup $\Gamma$, usually 
denoted by $m(\Gamma)$ and $e(\Gamma)$ respectively.
\medskip

The \textit{Ap\'{e}ry set} of $\Gamma$ with respect to a non-zero $a\in \Gamma$ is defined to be the set $\rm{Ap}(\Gamma,a)=\{s\in \Gamma\mid s-a\notin \Gamma\}$.
Each element $x \in \Gamma$ can be written as $x = \sum_{i=1}^{e} a_{i}s_{i}$
for some non-negative integers $s_{i}$. The vector $\mathbf{s}=(s_{1},\ldots,s_{e})$ is called a factorization of $x$, and the
set of all factorizations of $x$ is denoted by $F(x)$, which is a finite set.
\medskip

Let $\mid \mathbf{s}\mid =\sum_{i=1}^{e} s_{i}$ denote the total order of $\mathbf{s}$. Then the maximum integer $n$ which is the total order of a vector in $F(x)$ is called the order of $x$ and is denoted by $\mathrm{ord}_{\Gamma}(x)$.
A vector $\mathbf{s} \in F(x)$ with $\mid\mathbf{s}\mid = \mathrm{ord}_{\Gamma}(x)$, is called a maximal factorization of $x$ and $x = \sum_{i=1}^{e} a_{i}s_{i}$ is called a maximal expression of $x$. For a vector a of non-negative
integers, we set $x(\mathbf{s}) = \sum_{i=1}^{e} a_{i}s_{i}$. Given $0\neq x \in \Gamma$, the set of lengths of $x$ in $\Gamma$ is defined as
$$L(x) =\{\sum_{i=1}^{e}r_{i} | x = \sum_{i=1}^{e}r_{i}a_{i},r_i \geq 0\}$$
Given integers $ a_{1} < \cdots < a_{e}$; the map 
$\nu : k[x_{1}, \ldots, x_{e}]\longrightarrow k[t]$ defined as 
$\nu(x_{i}) = t^{a_{i}}$, $1\leq i\leq e$, defines a parametrization 
for an affine monomial curve; the ideal $\ker(\nu)=\mathfrak{p}$ is called the 
defining ideal of the monomial curve defined by the parametrization 
$\nu(x_{i}) = t^{a_{i}}$, $1\leq i\leq e$. The defining ideal 
$\mathfrak{p}$ is a graded ideal with respect to the weighted gradation 
and therefore any two minimal generating sets of $\mathfrak{p}$ have the 
same cardinality.
\medskip

Suppose $M=\Gamma\setminus\{0\}$ and for a positive integer $n$, we write 
$nM:=M+\cdots+M$ ($n$-copies). Let $r:=min\{r|(r+1)M=a_{1}+rM\}$, this $r$ is called the reduction number. Let $\mathfrak{m}$ be the maximal ideal of the 
ring $k[[t^{a_{1}},\ldots t^{a_{e}}]]$. Then $(n+1)M=a+nM$ for all $n\geq r$ 
if and only if $r=r_{(t^{a_{1}})}(\mathfrak{m})$.
\medskip

Let $\mathrm{Ap}(\Gamma,a_{1})=\{0,\omega_{1},\ldots,\omega_{a_{1}-1}\}$. Now for 
each $n\geq 1$, let us define $\mathrm{Ap}(nM)=\{\omega_{n,0},\ldots\omega_{n,a_{1}-1}\}$ 
inductively. We define $\omega_{1,0}=a_{1}$ and $\omega_{1,i}=\omega_{i}$, for $1\leq i\leq a_{1}-1$. 
Then $\mathrm{Ap}(M)=\{a_{1},\omega_{1},\ldots,\omega_{a_{1}-1}\}$. Now we define 
$\omega_{n+1,i}=\omega_{n,i}$, if $\omega_{n,i}\in (n+1)M $, and $\omega_{n+1,i}=\omega_{n,i}+a_{1}$, 
otherwise. We note that $\omega_{n+1,i}=\omega_{n,i}+a_{1} $ for all $0\leq i\leq a_{1}-1$ and $n\geq r_{(t^{a_{1}})}(\mathfrak{m})$. Then, the Ap\'{e}ry table $\mathrm{AT}(\Gamma,a_{1})$ of $\Gamma$ is a table 
of size $(r_{(t^{a_{1}})}(\mathfrak{m})+1)\times a_{1}$, whose $(0,t)$ entry is $\omega_{t}$, 
for $0\leq t\leq {a_{1}-1}$ (we take $\omega_{0}=0$), and the $(s,t)$ entry is $\omega_{st}$, 
for $1\leq s\leq r_{(t^{a_{1}})}(\mathfrak{m})$ and $ 0\leq t\leq {a_{1}-1}$.
\medskip

We take some definitions from \cite{cz2}. Let $W =\{a_{0},\ldots,a_{n}\}$ be a set of integers. We call it a \textit{ladder} if $a_{0}\leq\ldots\leq a_{n}$. Given a ladder, we say that a subset $L=\{a_{i},\ldots,a_{i+k}\}$, with $k\geq 1$, is a \textit{landing} of length $k$ if $a_{i-1}<a_{i}=\cdots=a_{i+k}<a_{i+k+1}$ (where $a_{-1}= -\infty$ and $a_{n+1}=\infty$). In this case, $s(L)=i$ and $e(L)=i+k$. A landing $L$ is said to be a \textit{true landing} if $s(L)\geq 1$. Given two landings $L$ and $L^{'}$, we set $L<L^{'}$ if $s(L)<s(L^{'})$. Let $p(W)+1$ be the number of landings and assume that $L_{0}<\cdots<L_{p(W)}$ are the distinct landings. Then we define the following numbers:
$s_{j}(W)=s(L_{j})$, $e_{j}(W)=e(L_{j})$, for each $0\leq j\leq p(W)$;
$c_{j}(W)=s_{j}(W)-e_{j-1}(W)$, for each $0\leq j\leq p(W)$.
\medskip

Suppose $\Gamma$ be a numerical semigroup minimally 
generated by $a_{1}<\cdots <a_{e}$ and $\mathfrak{m}_{\Gamma}$ be the maximal ideal of $k[[t^{a_{1}},\ldots t^{a_{e}}]]$. Let $r= r_{(t^{a_{1}})}(\mathfrak{m}_{\Gamma})$,  $M=\Gamma\setminus\{0\}$ and 
$\mathrm{Ap}(nM)=\{\omega_{n,0},\ldots\omega_{n,a_{1}-1}\}$ for $0\leq n \leq r$. For every $1\leq i\leq a_{1}-1$, consider the ladder of the values $W^{i}=\{\omega_{n,i}\}_{0\leq n\leq r}$ and define the following integers:
\begin{enumerate}[(i)]
\item $p_{i}=p(W^{i})$
\item $d_{i}=e_{p_{i}}(W^{i})$
\item $b_{j}^{i}=e_{j-1}(W^{i})$ and 
$c_{j}^{i}=c_{j}(W^{i})$, for $1\leq j\leq p_{i}$.
\end{enumerate}

\begin{theorem}\textbf{(Cortadellas, Zarzuela.)}\label{tangentcone} With the above notations, 
$$G_{\mathfrak{m}_{\Gamma}}\cong F_{\Gamma}\oplus\displaystyle\bigoplus_{i=1}^{a_{1}-1}\left(F_{\Gamma}(-d_{i})\displaystyle \bigoplus_{j=1}^{p_{i}}\dfrac{F_{\Gamma}}{(({t^{a_{1}})^{*})^{c_{j}^{i}}}F_{\Gamma}}(-b_{j}^{i})\right),$$
where $G_{\mathfrak{m}_{\Gamma}}$ is the tangent cone of $\Gamma$ and $F_{\Gamma}=F((t^{a_{1}}))$ is the fiber cone.
\end{theorem}

\proof See Theorem 2.3 in \cite{cz2}.\qed

\begin{theorem}Let $H_{F_{\Gamma} } (x)$ be the  Hilbert  series of $F_{\Gamma}$, then as stated in \cite{cz1},
$$F_{\Gamma } {\displaystyle \simeq }\displaystyle \bigoplus_{i=1}^{e} F (J )(-b_{i} )
\displaystyle \bigoplus_{j=1}^{f} (F (J )/a^{c_{j}} F (J ))(-d_{j}) $$
where  we  may  assume  $b_{1}\leq \cdots \leq  b_{e}$,  $d_{1}\leq \cdots \leq d_{f}$. In particular one immediately
has
$$ H_{F_{\Gamma}} (x)  =\dfrac{x^{b_{1}}+ \cdots + x^{b_{e}}       + (1 - x^{c_{1}})x^{d_{1}}+\cdots+ (1-x^{c_{f}})x^{d_{f}}}{
1-x} $$ 
\end{theorem}
As given in section $2$\cite{bjz}, we denote $H(n)=\# nM\setminus (n+1)M$ , the Hilbert function of $k[[\Gamma]]$

\section{Tangent cone of Bresinsky Curves}
Let $h\geq 2$ be an integer. Let $m_{0}=2h(2h-1)$, $m_{1}=(2h+1)(2h-1)$, $m_{2}=2h(2h+1)$, $m_{3}=2h(2h+1)+(2h-1)$. Bresinsky see \cite{bre} defined the curve $\Gamma_{h}=\langle m_{0},m_{1},m_{2},m_{3}\rangle$. 
\
\begin{theorem}\label{aperybre}
The Ap\'{e}ry set $Ap(\Gamma_{h},m_{0})$ is given as follows,
$$Ap(\Gamma_{h},m_{0})= A_{1}\cup A_{2}\cup A_{3}\cup A_{4}\cup A_{5}.$$
Where,
\begin{itemize}
\item $A_{1}=\{im_{1}|1\leq i\leq 2h-1\}$,
\item $A_{2}=\{im_{2}|1\leq i\leq 2h-2\}$,
\item  $A_{3}=\{im_{3}|1\leq i\leq 2h-2\}$,
\item $A_{4}=\{im_{1}+jm_{3}|1\leq i\leq 2h-2, 1\leq j\leq (2h-1)-i\}$,
\item $A_{5}=\{im_{2}+jm_{3}|1\leq i\leq 2h-3, 1\leq j\leq (2h-2)-i\}$
\end{itemize}

\end{theorem}
\proof  We have $m_{1}=m_{0}+(2h-1),\, m_{2}=m_{0}+4h,\, m_{3}= m_{0}+(6h-1)$ and we want to show that $A_{1}\subset Ap(\Gamma_{h},m_{0})$. At first, we show that $(2h-1)m_{1}-m_{0} \notin \Gamma_{h}$.
\medskip

Suppose $(2h-1)m_{1}-m_{0}= a_{0}m_{0}+a_{1}m_{1}+a_{2}m_{2}+a_{3}m_{3}$, hence
\begin{equation}\label{eq*}(2h-1)m_{1}=(a_{0}+1)m_{0}+a_{1}m_{1}+a_{2}m_{2}+a_{3}m_{3}. 
\end{equation} 
We claim that $a_{0}+a_{1}+a_{2}+a_{3} < 2h-1$.
\medskip

If $a_{0}+a_{1}+a_{2}+a_{3} \geq 2h-1 $ then $ (a_{0}+a_{1}+a_{2}+a_{3})m_{1} \geq (2h-1)m_{1}$. We have
\begin{eqnarray*}
(a_{0}+1)m_{0}-a_{0}m_{1} & = & (a_{0}+1)m_{0}-a_{0}(m_{0}+2h-1)\\
& = & m_{0}-a_{0}(2h-1)\\
& = & (2h-1)(2h-a_{0})
\end{eqnarray*}
If $a_{0} \geq 2h$, then $(a_{0}+1) \geq (2h+1)$ and
$$(a_{0}+1)m_{0} \geq (2h+1)2h(2h-1) \gneq (2h-1)^{2}(2h+1)=(2h-1)m_{1}$$
which implies $(2h-a_{0})>0$.
Therefore, $(a_{0}+1)m_{1} > a_{0}m_{1}$, which gives
\begin{align*}
(a_{0}+1)m_{0}+a_{1}m_{1}+a_{2}m_{2}+a_{3}m_{3} &> (a_{0}+a_{1}+a_{2}+a_{3})m_{1} \\
&\geq (2h-1)m_{1},
\end{align*}
 a contradiction. Therefore, \begin{equation}\label{eq**}
(a_{0}+1)+a_{1}+a_{2}+a_{3}  \leq  2h-1 
\end{equation}

Again from \ref{eq*}, $$(2h-1)[(2h-1-a_{1})+m_{0}]= (a_{0}+a_{1}+a_{2}+a_{3})m_{0}+a_{2}4h+a_{3}(6h-1)$$
which gives, $ 2h-1\mid a_{2}4h+a_{3}(6h-1)$. Therefore, $ 2h-1\mid 2(a_{2}+a_{3})$.
\medskip

Since $\gcd(2h-1,2)=1$, we have $ 2h-1 \mid (a_{2}+a_{3})$. Hence $a_{2}+a_{3}=2h-1 $ or $a_{2}+a_{3}=0$. If
$a_{2}+a_{3} = 2h-1$ , we get a contradiction RHS $>$ LHS in \ref{eq*}.
Therefore, $a_{2}+a_{3}=0$, i.e. $a_{2}=0=a_{3}$, which gives $(2h-1-a_{1})m_{1}=(a_{0}+1)m_{0}$ (from \ref{eq*}).
Hence $ 2h+1 \mid (a_{0}+1)m_{0}$. Since $\gcd(2h,2h+1)=1 $, and $\gcd(2h+1, 2h-1)=1$, we have $\gcd(m_{0},2h+1)=1$.
which implies $ 2h+1\mid (a_{0}+1)$ , gives a contradiction as $a_{0}+1 \leq 2h-1 $ by \ref{eq**}.
\medskip

If for any $1\leq i < 2h-1 $, $im_{1}-m_{0} \in \Gamma_{h}$, hence $ im_{1}-m_{0} + (2h-1-i)m_{1} \in \Gamma_{h}$, which implies $(2h-1)m_{1}-m_{0} \in  \Gamma_{h}$, a contradiction. Therefore, $A_{1} \subseteq A_{p}(\Gamma_{h},m_{0})$.
\smallskip

Similarly, we can show $A_{2}, A_{3}$ are also subsets of the Apery set.
\smallskip

Next we have, $A_{4}= \{ im_{1}+jm_{3} | 1\leq i \leq 2h-2, 1 \leq j \leq (2h-1)-i\}$ and it is enough to show that $ im_{1}+(2h-1-i)m_{3}-m_{0} \notin \Gamma_{h}$. Suppose 
\begin{equation}\label{eq***} 
im_{1}+ (2h-1-i)m_{3}= (a_{0}+1)m_{0}+a_{1}m_{1}+a_{2}m_{2}+a_{3}m_{3} 
\end{equation}
We claim that $a_{0}+a_{1}+a_{2}+a_{3} \leq 2h-1 $.\\
If $a_{0}+a_{1}+a_{2}+a_{3} \geq 2h $, we have
\begin{align*}
& im_{1}+ (2h-1-i)m_{3}-[(a_{0}+1)m_{0}+a_{1}m_{1}+a_{2}m_{2}+a_{3}m_{3}] \\
&=  i(m_{0}+2h-1)+(2h-1-i) (m_{0}+6h-1)\\ 
&-[(a_{0}+1)m_{0}+a_{1}(m_{0}+2h-1)+a_{2}(m_{0}+4h)+a_{3}(m_{0}+6h-1)] \\
&=  [(2h-1)-a_{0}-1-a_{1}-a_{2}-a_{3}]m_{0}\\
&+(i-a_{1})(2h-1)+(2h-1-i-a_{3})(6h-1)-a_{2}(4h)
\end{align*}
Maximum value of this expression is  $$-2m_{0}+(i-a_{1})(2h-1)+(2h-1-i-a_{3})(6h-1)-a_{2}4h < 0.$$ [ since  $a_{0}+a_{1}+a_{2}+a_{3} \geq 2h]$, 
Therefore LHS$ <$ RHS  of \ref{eq***} gives a contradiction. Again from \ref{eq***} we get 
\begin{align*}
& m_{0}(2h-1-a_{0}-1-a_{1}-a_{2}-a_{3})\\
&=(a_{1}-i)(2h-1)+(a_{3}-2h+1+i)(6h-1)+a_{2}4h,
\end{align*}
which implies $2h\mid (a_{1}-i)(2h-1)+(a_{3}-2h+1+i)(6h-1)+a_{2}4h$, hence $ 2h \mid -(a_{1}-i)-(a_{3}+1+i),$ i.e. $ 2h | a_{1}+a_{3}+1 $.Therefore, $a_{1}+a_{3}+1=0 $ or $ a_{1}+a_{3}+1 =2h $, as  $  a_{0}+a_{1}+a_{2}+a_{3}+1 \leq 2h$. But $a_{1}+a_{3}+1 \neq 0 $ as $ a_{1} \geq 0 \  or \ a_{3} \geq 0 $. Therefore, $a_{1}+a_{3}=2h-1$.
\medskip

The minimum value of RHS of \ref{eq***} is $(2h-1)m_{1}+m_{0}$. Now,
$$(2h-1)m_{1}+m_{0}-im_{1}-(2h-1-i)m_{3}= (2h-1)(4h-1-i)-(6h-1)>0$$ gives a contradiction. Therefore, $im_{1}+(2h-1-i)m_{3}-m_{0} \notin \Gamma_{h}$. Hence we can conclude that $A_{4} \subseteq A_{p}(\Gamma_{h},m_{0})$.
Similarly, we can prove that $A_{5} \subseteq A_{p}(\Gamma_{h},m_{0})$. \qed

\begin{theorem}\label{unique3}
Every element of $\mathrm{Ap}(\Gamma_{h},m_{0})$ has a unique expression.
\end{theorem}
\proof At first we will show that every element of the set $A_{1}$ in $\mathrm{Ap}(\Gamma_{h},m_{0})$, is uniquely expressed.  Since $m_{1}$ is the element of the generating set of the numerical semigroup $\Gamma_{h}$, it has a unique expression. Suppose $$(2h-1)m_{1}=c_{1}m_{1}+c_{2}m_{2}+c_{3}m_{3}.$$ If $c_{1}\neq 0$ then $(2h-1-c_{1})m_{1}=c_{2}m_{2}+c_{3}m_{3}$ and by induction $(2h-1-c_{1})m_{1}$ has a unique expression,a contradiction. 
\medskip

If $c_{1}=0$ then $(2h-1)m_{1}=c_{2}m_{2}+c_{3}m_{3}$ taking modulo $2h-1$ we get,
$$0\equiv 2(c_{2}+c_{3})\mod (2h-1).$$ Hence, $2(c_{2}+c_{3})=k(2h-1)$. By the equation $(2h-1-c_{2}-c_{3})m_{1}=c_{2}(2h+1)+c_{3}4h$, we get $c_{2}+c_{3}<(2h-1)$. The LHS $2(c_{2}+c_{3})<2(2h-1)$, and RHS $k_{1}(2h-1)\geq 2(2h-1)$, which is a contradiction, therefore the elements of $A_{1}$ are uniquely expressed.
\medskip

 Since $m_{2}$ is the element of the generating set of the numerical semigroup $\Gamma_{h}$, so it has a unique expression. Suppose $(2h-2)m_{2}=c_{1}m_{1}+c_{2}m_{2}+c_{3}m_{3}$. If $c_{2}\neq 0$ then $(2h-2-c_{2})m_{2}=c_{1}m_{1}+c_{3}m_{3}$ and by induction $(2h-2-c_{2})m_{2}$ has a unique expression,a contradiction. If $c_{2}=0$ then $(2h-2)m_{2}=c_{1}m_{1}+c_{3}m_{3}$ taking modulo $2h-1$ we get,$$-2\equiv 2c_{3}\mod (2h-1).$$ Hence $2c_{3}+2=k_{2}(2h-1)$. If $k_{2}=1$, then $c_{3}=\dfrac{2h-3}{2}$. Therefore, we get $k_{2}\geq 2$.
\medskip

 Again we have $$(2h-2-c_{1}-c_{3})m_{1}+(2h-2)(2h+1)=c_{3}4h.$$ If $c_{1}+c_{3}>(2h-2)$ then $(2h-2-c_{1}-c_{3})m_{1}+(2h-2)(2h+1)<-(2h+1)$ and RHS is $c_{3}4h\geq 0$. Therefore,  $c_{1}+c_{3}\leq (2h-2)$. If $c_{1}+c_{3}=(2h-2)$ and $c_{1}=0$, then $(2h+1)=4h$, which is a contradiction. Hence $c_{3}< (2h-2)$. In that case, LHS is $2c_{3}+2< 4h-2$ and RHS is $k_{2}(2h-1)\geq 4h-2$, which is a contradiction.
Similarly, we can show that the elements of $A_{3}$ are uniquely expressed. 
\medskip

Let $im_{1}+(2h-1-i)m_{3}=c_{1}m_{1}+c_{2}m_{2}+c_{3}m_{3}$, Then, we get
$(2h-1-c_{1}-c_{2}-c_{3})m_{1}+(2h-1-i)4h=c_{2}(2h+1)+c_{3}4h$.
\medskip

If $c_{1}+c_{2}+c_{3}>(2h-1)$, then LHS is $(2h-1-c_{1}-c_{2}-c_{3})m_{1}+(2h-1-i)4h< -8h+2$. Since $h\geq 2$ therefore LHS$<-14$ and RHS $c_{2}(2h+1)+c_{3}4h>0$, which is a contradiction.
\medskip 

If $c_{1}+c_{2}+c_{3}=2h$, then $-m_{1}+(2h-1-i-c_{3})4h=c_{2}(2h+1)$.
Since RHS is $c_{2}(2h+1)>0$ then $(2h-1-i-c_{3})>0$ and $(2h+1)|(2h-1-i-c_{3})$. Since $(2h+1)>(2h-1-i-c_{3})$, which is a contradiction. 
\medskip

If $c_{1}+c_{2}+c_{3}=(2h-1)$, then $ (2h-1-i-c_{3})4h=c_{2}(2h+1)$. Hence $c_{2}=k.2h$. As $c_{1}+c_{2}+c_{3}=(2h-1)$, therefore $c_{2}=0$ hence $im_{1}+(2h-1-i)m_{3}=c_{1}m_{1}+c_{3}m_{3}$ and $c_{1}+c_{3}=(2h-1)$. We already have that expression. 
\medskip

If $c_{1}+c_{2}+c_{3}<(2h-1)$, then $(2h-1)-c_{3}>c_{1}+c_{2}\geq 0$ and $(2h-1)-c_{3}-i>-i$. We know that $2h+1|(2h-1)-c_{3}-i$. Since $1\leq i \leq (2h-2)$, therefore $(2h-1)-c_{3}-i=0$. Substituting the value of $c_{3}$ in
 $(2h-1-c_{1}-c_{2}-c_{3})m_{1}+(2h-1-i)4h=c_{2}(2h+1)+c_{3}4h$, we get $(2h-1-c_{1}-c_{2}-c_{3})m_{1}=c_{2}(2h+1)$. Therefore, $2h-1|c_{2}$, let $c_{2}=k(2h-1)$. Since $c_{1}+c_{2}+c_{3}<(2h-1)$, therefore, $k=0$. Hence $(2h-1)=c_{1}-c_{2}-c_{3}$, which is a contradiction.
 \medskip
 
 Similarly, we can show that elements of $A_{5}$ are uniquely expressed. \qed

\begin{theorem}\label{aperytable3}
The following statements hold for all $k\geq 0$.
\begin{enumerate}
\item $\mathrm{ord}(im_{1}+km_{0})=i+k$, $1 \leq i \leq 2h-1$
\item $\mathrm{ord}(im_{2}+km_{0})=i+k$, $1 \leq i \leq 2h-2$
\item $\mathrm{ord}(im_{3}+km_{0})=i+k$, $1 \leq i \leq 2h-2$
\item $\mathrm{ord}(im_{1}+jm_{3}+km_{0})=i+j+k$, $1 \leq i \leq 2h-2$, $1 \leq j \leq (2h-1)-i$
\item $\mathrm{ord}(im_{2}+jm_{2}+km_{0})=i+j+k$, $1 \leq i \leq 2h-3$ , $1 \leq j \leq (2h-2)-i$
\end{enumerate}
\end{theorem}

\proof In $(1)$, we want to show $\mathrm{ord}(im_{1}+km_{0})=i+k$, for all $k \geq 0$ , $1 \leq i \leq 2h-1$.
Suppose that we prove $\mathrm{ord}((2h-1)m_{1}+km_{0})=2h-1+k$ for all  $k \geq 0$ and if $\mathrm{ord}(im_{1}+km_{0})>i+k$ for some  $1 \leq i \leq 2h-1$ , then $(2h-1)m_{1}+km_{0}=(2h-1-i)m_{1}+im_{1}+km_{0}$.
i.e. $ord((2h-1)m_{1}+km_{0})> 2h-1+k$ that is a contradiction.
Therefore, it is enough to show, $\mathrm{ord}((2h-1)m_{1}+km_{0})=2h-1+k$ for all $k \geq 0$.
\medskip

Suppose, $(2h-1)m_{1}+km_{0}= a_{0}m_{0}+a_{1}m_{1}+a_{2}m_{2}+a_{3}m_{3}$.
If $a_{0} \geq k $, $(2h-1)m_{1}=(a_{0}-k)m_{0}+a_{1}m_{1}+a_{2}m_{2}+a_{3}m_{3}$. But $(2h-1)m_{1} \in A_{p}(\Gamma_{h},m_{0})$, has a unique expression by the theorem \ref{unique3}. Therefore, $a_{1}=2h-1$ and $a_{0}=k$, $a_{2}=0$, $a_{3}=0$ and we are done.
\medskip

If $0<a_{0}<k$ then, we have $$(2h-1)m_{1}+(k-a_{0})m_{0}= a_{1}m_{1}+a_{2}m_{2}+a_{3}m_{3}$$ and $k-a_{0} < k$, by induction hypothesis we are done.
\medskip

So the only case is, $a_{0}=0$. We have $$(2h-1)m_{1}+km_{0}=a_{1}m_{1}+a_{2}m_{2}+a_{3}m_{3}.$$

If $a_{1}\geq (2h-1)$, then $km_{0}=(a_{1}-(2h-1))m_{1}+a_{2}m_{2}+a_{3}m_{3}$. Since $m_{0}<m_{1}<m_{2}<m_{3}$, we have
$(a_{1}-(2h-1))m_{0} \leq (a_{1}-(2h-1))m_{1}$, $a_{2}m_{0}\leq a_{2}m_{2}$ and 
$a_{3}m_{0} \leq a_{3}m_{3}$. Therefore, by adding these, we get
$(a_{1}-(2h-1)+a_{2}+a_{3})m_{0} \leq km_{0}$. which implies that $ a_{1}-(2h-1)+a_{2}+a_{3} \leq k$, hence $ a_{1}+a_{2}+a_{3} \leq (2h-1)+k $ and we are done.
\medskip

If $0 \leq a_{1} < (2h-1)$ then, $((2h-1)-a_{1})m_{1}+km_{0}=a_{2}m_{2}+a_{3}m_{3}$. We have, $((2h-1)-a_{1})m_{1} \leq ((2h-1)-a_{1})m_{2}$ and $ km_{0} \leq km_{2}$. Therefore, by adding these two equations, we get  $$((2h-1)-a_{1})m_{1}+km_{0} \leq ((2h-1)-a_{1})m_{2}+km_{2}.$$
\medskip

We have
\begin{align*}
(a_{2}+a_{3})m_{2}& \leq a_{2}m_{2}+a_{3}m_{3} \\
&\leq ((2h-1)-a_{1}+k)m_{2}.
\end{align*}
Which implies, $a_{2}+a_{3} \leq (2h-1)-a_{1}+k $, i.e., $a_{1}+a_{2}+a_{3} \leq 2h-1+k $, and we are done.The proof of the statements $(2), (3), (4), \text{ and } (5) $ are similar as above. \qed

\begin{theorem}\label{atbre}
 The Ap\'{e}ry table of the Bresinsky curve w.r.t. $m_{0}$ is a matrix of order $2h\times m_{0}$ and the table $\mathrm{AT}(\Gamma_{h},m_{0})$ is given as follows,
$$ \mathrm{AT}(\Gamma_{h},m_{0})=\begin{bmatrix}
  \leftidx{^h}{\mathcal{T}_{0}}  &\leftidx{^h}{\mathcal{T}_{1}} &\leftidx{^h}{\mathcal{T}_{2}} &\leftidx{^h}{\mathcal{T}_{3}} &\leftidx{^h}{\mathcal{T}_{4}} &\leftidx{^h}{\mathcal{T}_{5}}
 \end{bmatrix},$$  where 
 \begin{itemize}
     \item $\leftidx{^h}{\mathcal{T}_{0}}=(t^{(0)}_{i1})_{2h\times 1}$ and $t^{(0)}_{i1}=(i-1)m_{0}$ for $1\leq i\leq 2h$
     \item $\leftidx{^h}{\mathcal{T}_{1}}=(t^{(1)}_{ij})_{2h\times (2h-1)}$ and 
   \begin{align*}
       t^{(1)}_{ij}&=jm_{1} \quad \mathrm{for}\,1\leq i\leq j+1, \,\, 1\leq j\leq 2h-1\\
      &=jm_{1}+(i-j-1)m_{0} \quad\mathrm{for}\,j+1< i\leq 2h, 1\leq j\leq 2h-1.
   \end{align*}  
   
     \item For $l=2,3$, $\leftidx{^h}{\mathcal{T}_{l}}=(t^{(l)}_{ij})_{2h\times (2h-2)}$ and 
   \begin{align*}
       t^{(l)}_{ij}&=jm_{l} \quad \mathrm{for}\,1\leq i\leq j+1, \,\, 1\leq j\leq 2h-2\\
      &=jm_{l}+(i-j-1)m_{0} \quad\mathrm{for}\, j+1< i\leq 2h, 1\leq j\leq 2h-2.
   \end{align*} 
    \item $\leftidx{^h}{\mathcal{T}_{4}}=\begin{bmatrix}
  \leftidx{^h}{T^{(4)}_{1}} &\leftidx{^h}{T^{(4)}_{2}} &\cdots  &\leftidx{^h}{T^{(4)}_{2h-2}} 
 \end{bmatrix}$ is a matrix of order $2h\times (h-1)(2h-1)$ and  $\leftidx{^h}{T^{(4)}_{i}}=(t^{(i4)}_{rs})_{2h\times (2h-1-i)}$ where,
   \begin{align*}
       t^{(i4)}_{rs}&=im_{1}+sm_{3}; \,1\leq r\leq i+s+1, \,\, 1\leq s\leq 2h-1-i\\
      &=im_{1}+sm_{3}+(r-i-s-1)m_{0};\\
      &\mathrm{for}\,\, r+s+1< i\leq 2h, 1\leq s\leq 2h-1-i.
\end{align*}  
\item $\leftidx{^h}{\mathcal{T}_{5}}=\begin{bmatrix}
  \leftidx{^h}{T^{(5)}_{1}} &\leftidx{^h}{T^{(5)}_{2}} &\cdots  &\leftidx{^h}{T^{(5)}_{2h-3}} 
 \end{bmatrix}$ is a matrix of order $2h\times (h-1)(2h-3)$ and $\leftidx{^h}{T^{(5)}_{i}}=(t^{(i5)}_{rs})_{2h\times (2h-2-i)}$ where,
   \begin{align*}
       t^{(i5)}_{rs}&=im_{2}+sm_{3}; \,1\leq r\leq i+s+1, \,\, 1\leq s\leq 2h-2-i\\
      &=im_{2}+sm_{3}+(r-i-s-1)m_{0};\\
      &\mathrm{for}\,\, r+s+1< i\leq 2h, 1\leq s\leq 2h-2-i.
\end{align*}  
 \end{itemize}
 \end{theorem}
 
\proof The proof follows from \ref{aperytable3}.\qed

\begin{example}
We take $h=3$, then $\Gamma_{3}=\langle 30,35,42,47\rangle$ and the Ap\'{e}ry table of the Bresinsky curve w.r.t. $m_{0}=30$ is a matrix of order $6\times 30$. Here $ \mathrm{AT}(\Gamma_{3},30)=\begin{bmatrix}
  \leftidx{^3}{\mathcal{T}_{0}}  &\leftidx{^3}{\mathcal{T}_{1}} &\leftidx{^3}{\mathcal{T}_{2}} &\leftidx{^3}{\mathcal{T}_{3}} &\leftidx{^3}{\mathcal{T}_{4}} &\leftidx{^3}{\mathcal{T}_{5}}
 \end{bmatrix}$ and
\begin{itemize}
    \item $\leftidx{^3}{\mathcal{T}_{0}}=\begin{bmatrix}
        0\\
        30\\
        60\\
        90\\
        120\\
        150
    \end{bmatrix}$ ,$\bullet\,\,\leftidx{^3}{\mathcal{T}_{1}}=\begin{bmatrix}
        35 & 70& 105 & 140& 175\\
        35 & 70& 105 & 140& 175\\
        65 & 70& 105 & 140& 175\\
        95 & 100& 105 & 140& 175\\
       125 & 130& 135 & 140& 175\\
        155 & 160& 165 & 175& 175\\
    \end{bmatrix}$ 
    \item $\leftidx{^3}{\mathcal{T}_{2}}=\begin{bmatrix}
        42 & 84& 126 & 168\\
        42 & 84& 126 & 168\\
        72 & 84& 126 & 168\\
        102 & 114& 126 & 168\\
       132 & 144& 156 & 168\\
        162 & 174& 186 & 198\\
    \end{bmatrix}$ $\bullet\,\,\leftidx{^3}{\mathcal{T}_{3}}=\begin{bmatrix}
        47 & 94& 141 & 188\\
         47 & 94& 141 & 188\\
         77 & 94& 141 & 188\\
        107 & 124& 141 & 188\\
       137 & 154& 171 & 188\\
         167 & 184& 201 & 218\\
    \end{bmatrix}$ 
    \item $\leftidx{^3}{\mathcal{T}_{4}}=\begin{bmatrix}
  \leftidx{^3}{T^{(4)}_{1}} &\leftidx{^3}{T^{(4)}_{2}} &\leftidx{^3}{T^{(4)}_{3}}  &\leftidx{^3}{T^{(4)}_{4}} 
 \end{bmatrix}$ is a matrix of order $6\times 8$ and \\
 $\bullet\,\leftidx{^3}{T^{(4)}_{1}}=\begin{bmatrix}
        82 & 129 & 176  & 223\\
       82 & 129 & 176  & 223\\
       82 & 129 & 176  & 223\\
       112 & 129 & 176  & 223\\
       142 & 159 & 176  & 223\\
      172  & 189 & 206  & 223\\
    \end{bmatrix}$,  $\bullet\,\,\leftidx{^3}{T^{(4)}_{2}}=\begin{bmatrix}
       117 &164 & 211\\
        117 &164 & 211\\
         117 &164 & 211\\
          117 &164 & 211\\
           147 &164 & 211\\
            177 &194 & 211\\
    \end{bmatrix}$,  $\bullet\,\,\leftidx{^3}{T^{(4)}_{3}}=\begin{bmatrix}
       152 & 199\\
       152 & 199\\
       152 & 199\\
       152 & 199\\
       152 & 199\\
       182 & 199\\
    \end{bmatrix}$ $\bullet\,\,\leftidx{^3}{T^{(4)}_{4}}=\begin{bmatrix}
        187\\
        187\\
        187\\
        187\\
        187\\
        187
    \end{bmatrix}$.
     \item $\leftidx{^3}{\mathcal{T}_{5}}=\begin{bmatrix}
  \leftidx{^3}{T^{(5)}_{1}} &\leftidx{^3}{T^{(5)}_{2}} &\leftidx{^3}{T^{(5)}_{3}}
 \end{bmatrix}$ is a matrix of order $6\times 6$ and \\
 $\bullet\,\,\leftidx{^3}{T^{(5)}_{1}}=\begin{bmatrix}
        89 & 136 & 183 \\
      89 & 136 & 183 \\
       89 & 136 & 183 \\
       119 & 136 & 183 \\
       149 & 166 & 183 \\
      179 & 196 & 213 \\
    \end{bmatrix}$,    $\bullet\,\,\leftidx{^3}{T^{(5)}_{2}}=\begin{bmatrix}
       131 & 178\\
       131 & 178\\
       131 & 178\\
       131 & 178\\
       161 & 178\\
       191 & 208\\
    \end{bmatrix}$ \\ $\bullet\,\,\leftidx{^3}{T^{(5)}_{3}}=\begin{bmatrix}
        173\\
        173\\
        173\\
        173\\
        173\\
        203
    \end{bmatrix}$.
\end{itemize}
\end{example}
\medskip

\begin{corollary}\label{order3}
 Let $I=\langle t^{m_{0}}\rangle$. The tangent cone $G_{\mathfrak{m}}(\Gamma_{h})$ of $\Gamma_{h}$ is a free $F(I)$-module. Moreover, 
$$ G_{\mathfrak{m}}(\Gamma_{h})= \displaystyle\bigoplus_{k=1}^{2h-2}(F(I)(-k))^{2k+1}\bigoplus(F(I)(2h-1))^{2h-1}.$$
\end{corollary}
\proof The proof follows from the theorems \ref{tangentcone}, \ref{unique3}, \ref{aperytable3} and lemma \ref{order3}.\qed

\begin{corollary}\label{CMBre}
The tangent cone $G_{\mathfrak{m}}(\Gamma_{h})$ is Cohen-Macaulay.
\end{corollary}
\proof It is easily followed from the fact that $G_{\mathfrak{m}}(\Gamma_{h})$ is a free $F(I)$-module (see section 4 in \cite{cz2}). \qed

\begin{corollary} Let $HG_{\mathfrak{m}}(\Gamma_{h})(x)$ be the Hilbert series of $G_{\mathfrak{m}}$. Then $$HG_{\mathfrak{m}}(\Gamma_{h})(x)=\displaystyle \left(\sum_{k=1}^{2h-2} (2k+1)x^{k}+(2h-1)x^{2h-1}\right)/(1-x)$$.
\end{corollary}
\proof The proof follows from Corollary \ref{CMBre}.\qed

\begin{remark}
Cohen-Macaulayness of the tangent cone Bresinsky curves has been already studied in \cite{hd}. But here we study the Ap\'{e}ry table and we give an explicit structure of the tangent cone of Bresinsky curves.
\end{remark}

\section{Tangent cone of Arslan curves}
Let $m\geq 2$ and $n_{1}=m(m+1),n_{2}=m(m+1)+1, n_{3}=(m+1)^{2}, n_{4}=(m+1)^{2}+1$. Arslan in \cite{ars} defined the following curves $\mathfrak{S}_{m}=\langle n_{1},n_{2},n_{3},n_{4}\rangle$.
\begin{theorem}\label{unique4}
 The Ap\'{e}ry set $\mathrm{Ap}(\mathfrak{S}_{m},n_{1})$ is given as follows
$$\mathrm{Ap}(\mathfrak{S}_{m},n_{1})=\mathfrak{A}_{1}\cup \mathfrak{A}_{2}\cup \mathfrak{A}_{3}\cup\mathfrak{A}_{4}\cup \mathfrak{A}_{5}$$ where 
\begin{itemize}
\item $\mathfrak{A}_{1}= \{in_{2}|1 \leq i \leq m \}$
\item $\mathfrak{A}_{2}= \{in_{3}|1 \leq i \leq m-1 \}$
\item $\mathfrak{A}_{3}= \{in_{4}|1 \leq i \leq m-1 \}$
\item $\mathfrak{A}_{4}= \{in_{2}+jn_{4}|1 \leq i \leq m-1,1 \leq j \leq m-i \}$
\item $\mathfrak{A}_{5}= \{in_{3}+jn_{4}|1 \leq i \leq m-2,1 \leq j \leq (m-1)-i \}$
\end{itemize}
\end{theorem}

\proof The proof is similar as the Theorem \ref{aperybre}. For example, if we want to show $\mathfrak{A}_{1}\subseteq A_{p}(\mathfrak{S}_{m},n_{1})$.Then it is enough to show $mn_{2}-n_{1} \notin \mathfrak{S}_{m}$.
Suppose $mn_{2}-n_{1} =a_{1}n_{1}+a_{2}n_{2}+a_{3}n_{3}+a_{4}n_{4}$. Then \begin{equation}\label{equ*}
[m-(a_{1}+1)-a_{2}-a_{3}-a_{4}]n_{1}+m=a_{2}+a_{3}(m+1)+a_{4}(m+2) 
\end{equation}
Since RHS of the equation \ref{equ*} is $\geq 0$. We have $m-(a_{1}+1)-a_{2}-a_{3}-a_{4} \geq 0 $. If $m-(a_{1}+1)-a_{2}-a_{3}-a_{4} = 0 $,
then from \ref{equ*}, we get $m=a_{2}+a_{3}(m+1)+a_{1}(m+2)$. Which implies
$ a_{2}=m, a_{3}=0, a_{4}=0 $ and substituting this values in \ref{equ*} we get
$\ mn_{2}=(a_{1}+1)n_{1}+mn_{2}$, a contradiction as $a_{1} \geq 0$. Therefore, we have 
\begin{equation}\label{equ**}
(a_{1}+1)+a_{2}+a_{3}+a_{4} < m 
\end{equation}
Again from \ref{equ*}, we get $(m+1)\mid a_{2}+a_{3}(m+1)+a_{4}(m+2)-m$
i.e. $m+1 \mid (a_{2}+a_{4})-m$, but $0<a_{2}+a_{4}<m$ gives a contradiction.
By a similar method as in \ref{aperybre},
$\mathfrak{A}_{i} \in \mathrm{Ap}(\mathfrak{S}_{m},n_{1})$, $2 \leq i \leq 5$ \qed

\begin{theorem}\label{uniquearslan}
Each element of $\mathrm{Ap}(\mathfrak{S}_{m},n_{1})$ has a unique expression.
\end{theorem}

\proof At first we will show that every element of the set $\mathfrak{A}_{1}$ in $Ap(\mathfrak{S}_{m},n_{1})$ is uniquely expressed.  Since $n_{2}$ is the element of the generating set of the numerical semigroup $\mathfrak{S}_{m}$, so it has a unique expression. Suppose $mn_{2}=c_{2}n_{1}+c_{3}n_{3}+c_{4}n_{4}$. If $c_{2}\neq 0$ then $(m-c_{2})n_{2}=c_{3}n_{3}+c_{4}n_{4}$ and by induction $(m-c_{2})n_{2}$ has a unique expression,a contradiction. If $c_{2}=0$, then $m n_{2}=c_{3}(n_{2}+m)+c_{4}(n_{2}+m+1)$ which gives, $(m-c_{3}-c_{4})n_{2}=c_{3}m+c_{4}(m+1)$.
\smallskip

If $c_{3}+c_{4}>m$, then L.H.S., $(m-c_{3}-c_{4})n_{2}<-{n_{2}}$ and R.H.S., $c_{3}m+c_{4}(m+1)>0$ which is a contradiction.
\smallskip

If $c_{3}+c_{4}=m$, then L.H.S., $(m-c_{3}-c_{4})n_{2}=0$ and R.H.S., $c_{3}m+c_{4}(m+1)>0$, which is a contradiction.
\smallskip

If $c_{3}+c_{4}<m$, then L.H.S. $(m-c_{3}-c_{4})n_{2}\geq n_{2}=m(m+1)+1$ R.H.S $c_{3}m+c_{4}(m+1)<m^{2}+m$, which is a contradiction. Similarly, we can prove that the elements of $\mathfrak{A}_{2}$ and $\mathfrak{A}_{3}$ are uniquely expressed.
\smallskip

Let $in_{2}+(m-i)n_{4}=c_{2}n_{2}+c_{3}n_{3}+c_{4}n_{4}$, we get $(m-c_{2}-c_{3}-c_{4})n_{2}+(m-i)(m+1)=c_{3}m+c_{4}(m+1)$
\smallskip

If $c_{2}+c_{3}+c_{4}>m$, then L.H.S, $(m-c_{2}-c_{3}-c_{4})n_{2}+(m-i)(m+1)<-i(m+1)-1<0$ and R.H.S, $c_{3}m+c_{4}(m+1)>0$, which is a contradiction.
\smallskip

If $c_{2}+c_{3}+c_{4}=m$, then $(m-i-c_{4})(m+1)=c_{3}m$, $c_{3}=k(m+1)$. Since $c_{2}+c_{3}+c_{4}=m$, therefore $c_{3}=0$. Hence $c_{2}+c_{4}=m$ and $in_{2}+(m-i)n_{4}=c_{2}n_{2}+c_{4}n_{4}$, we already have that expression. 
\smallskip

If $c_{2}+c_{3}+c_{4}<m$, then $(m-i-c_{4})(m+1)\geq (m+1)^{2}+1$ and R.H.S., $c_{3}m+c_{4}(m+1)<m(m+1)$, which is a contradiction. 
\smallskip

Similarly, we can show that elements of $\mathfrak{A}_{5}$ are uniquely expressed.
\qed

\begin{theorem}\label{aperytable4}
The following statements holds; for all $k \geq 0$
\begin{enumerate}
\item $\mathrm{ord}(in_{2}+kn_{1})=i+k$ , $1 \leq i \leq m$
\item $\mathrm{ord}(in_{3}+kn_{1})=i+k$ , $1 \leq i \leq m-1$
\item $\mathrm{ord}(in_{4}+kn_{1})=i+k$ , $1 \leq i \leq m-1$
\item $\mathrm{ord}(in_{2}+jn_{4}+kn_{1})=i+j+k$ , $1 \leq i \leq m-1$, $1 \leq j \leq m-i$
\item $\mathrm{ord}(in_{3}+jn_{4}+kn_{1})=i+j+k$ , $1 \leq i \leq m-2$, $1 \leq j \leq (m-1)-i$
\end{enumerate}
\end{theorem}
\proof The proof is similar as in \ref{aperytable3}. \qed
\begin{theorem}\label{atars}
 The Ap\'{e}ry table of the Arslan curve w.r.t. $n_{1}$ is a matrix of order $(m+1)\times n_{1}$ and the table $\mathrm{AT}(\mathfrak{S}_{m},n_{1})$ is given as follows,
$$ \mathrm{AT}(\mathfrak{S}_{m},n_{1})=\begin{bmatrix}
  \leftidx{^m}{\mathcal{A}_{0}}  &\leftidx{^m}{\mathcal{A}_{1}} &\leftidx{^m}{\mathcal{A}_{2}} &\leftidx{^m}{\mathcal{A}_{3}} &\leftidx{^m}{\mathcal{A}_{4}} &\leftidx{^m}{\mathcal{A}_{5}}
 \end{bmatrix},$$  where 
 \begin{itemize}
     \item $\leftidx{^m}{\mathcal{A}_{0}}=(a^{(0)}_{i1})_{(m+1)\times 1}$ and $a^{(0)}_{i1}=(i-1)n_{1}$ for $1\leq i\leq (m+1)$
     \item   $\leftidx{^m}{\mathcal{A}_{1}}=(a^{(1)}_{ij})_{(m+1)\times m}$ and 
   \begin{align*}
       a^{(1)}_{ij}&=jn_{2} \quad \mathrm{for}\,1\leq i\leq j+1, \,\, 1\leq j\leq m\\
      &=jn_{2}+(i-j-1)n_{1} \quad\mathrm{for}\,j+1< i\leq m+1, 1\leq j\leq m.
   \end{align*} 
     \item For $l=2,3$,  $\leftidx{^m}{\mathcal{A}_{l}}=(a^{(l)}_{ij})_{(m+1)\times (m-1)}$ and 
   \begin{align*}
       a^{(l)}_{ij}&=jn_{l+1} \quad \mathrm{for}\,1\leq i\leq j+1, \,\, 1\leq j\leq m-1\\
      &=jn_{l+1}+(i-j-1)n_{1} \quad\mathrm{for}\,j+1< i\leq m+1, 1\leq j\leq m-1.
   \end{align*}  
\item $\leftidx{^m}{\mathcal{A}_{4}}=\begin{bmatrix}
  \leftidx{^m}{A^{(4)}_{1}} &\leftidx{^m}{A^{(4)}_{2}} &\cdots  &\leftidx{^m}{A^{(4)}_{m-1}} 
 \end{bmatrix}$ is a matrix of order $(m+1)\times \frac{m(m-1)}{2}$ and  $\leftidx{^m}{A^{(4)}_{i}}=(a^{(i4)}_{rs})_{(m+1)\times (m-i)}$ where,
\begin{align*}
a^{(i4)}_{rs}&=in_{2}+sn_{4};\,1\leq r\leq i+s+1,\, 1\leq s\leq m-i\\
&=in_{2}+sm_{4}+(r-i-s-1)n_{1}\\ 
&\mathrm{for}\,\, r+s+1< i\leq m+1, 1\leq s\leq m-i.
\end{align*}  
   
\item $\leftidx{^m}{\mathcal{A}_{5}}=\begin{bmatrix}
  \leftidx{^m}{A^{(5)}_{1}} &\leftidx{^m}{A^{(5)}_{2}} &\cdots  &\leftidx{^m}{A^{(5)}_{m-2}} 
 \end{bmatrix}$ is a matrix of order $(m+1)\times \frac{(m-2)(m-1)}{2}$ and  $\leftidx{^m}{A^{(5)}_{i}}=(a^{(i5)}_{rs})_{(m+1)\times (m-1-i)}$ where,
\begin{align*}
a^{(i5)}_{rs}&=in_{3}+sn_{4}; \,1\leq r\leq i+s+1, \, 1\leq s\leq m-1-i\\
&=in_{3}+sm_{4}+(r-i-s-1)n_{1};\\
&\mathrm{for}\,\, r+s+1< i\leq m+1, 1\leq s\leq m-1-i.
\end{align*}  
 \end{itemize}
\end{theorem}

\proof The proof follows from the theorem \ref{aperytable4}.\qed

\begin{example}
We take $m=4$, then $\mathfrak{S}_{4}=\langle 20,21,25,26\rangle$ and the Ap\'{e}ry table of the Arslan curve w.r.t. $n_{1}=20$ is a matrix of order $5\times 20$. Here $ \mathrm{AT}(\mathfrak{S}_{4},20)=\begin{bmatrix}
    \leftidx{^4}{\mathcal{A}_{0}}  &\leftidx{^4}{\mathcal{A}_{1}} &\leftidx{^4}{\mathcal{A}_{2}} &\leftidx{^4}{\mathcal{A}_{3}} &\leftidx{^4}{\mathcal{A}_{4}} &\leftidx{^4}{\mathcal{A}_{5}}
 \end{bmatrix}$ and
\begin{itemize}
    \item $\leftidx{^4}{\mathcal{A}_{0}}=\begin{bmatrix}
        0\\
        20\\
        40\\
        60\\
        80
    \end{bmatrix}$ ,$\bullet\,\,\leftidx{^4}{\mathcal{A}_{1}}=\begin{bmatrix}
        21 & 42& 63 & 84\\
        21 & 42& 63 & 84\\
        41 & 42& 63 & 84\\
        61 & 62& 63 & 84\\
      81 & 82& 83 & 84\\
      
    \end{bmatrix}$ 
    \item $\leftidx{^4}{\mathcal{A}_{2}}=\begin{bmatrix}
        25 & 50& 75\\
        25 & 50& 75\\
        45 & 50& 75\\
       65 & 70& 75\\
        85 & 90& 95
    \end{bmatrix}$ $\bullet\,\,\leftidx{^4}{\mathcal{A}_{3}}=\begin{bmatrix}
        26 & 52& 78\\
        26 & 52& 78\\
        46 & 52& 78\\
       66 & 72& 78\\
         86 & 92& 98\\
    \end{bmatrix}$ 
    \item $\leftidx{^4}{\mathcal{A}_{4}}=\begin{bmatrix}
   \leftidx{^4}{A^{(4)}_{1}} &\leftidx{^4}{A^{(4)}_{1}} &\leftidx{^4}{A^{(4)}_{1}}   
 \end{bmatrix}$ is a matrix of order $5\times 6$ and \\
 $\bullet\,\,\leftidx{^4}{A^{(4)}_{1}}=\begin{bmatrix}
        47 & 73 & 99\\
      47 & 73 & 99\\
       47 & 73 & 99\\
      67 & 73 & 99\\
   87 & 93 & 99
    \end{bmatrix}$,  $\bullet\,\,\leftidx{^4}{A^{(4)}_{2}}=\begin{bmatrix}
       68 &94 \\
        68 &94 \\
           68 &94 \\
           68 &94 \\
             88 &94 
    \end{bmatrix}$,  $\bullet\,\,\leftidx{^4}{A^{(4)}_{3}}=\begin{bmatrix}
       89\\
       89\\
      89\\
       89\\
       89
    \end{bmatrix}$ 
     \item $\mathcal{A}_{5}=\begin{bmatrix}
   A^{(5)}_{1} &A^{(5)}_{2}  
 \end{bmatrix}$ is a matrix of order $5\times 3$ and \\
   $\bullet\,\,\leftidx{^4}{A^{(5)}_{1}}=\begin{bmatrix}
       51 & 77\\
       51 & 77\\
      51 & 77\\
       71 & 77\\
       91 & 97\\
       
    \end{bmatrix}$, $\bullet\,\,\leftidx{^4}{A^{(5)}_{2}}=\begin{bmatrix}
        76\\
        76\\
        76\\
        76\\
        96\\
     
    \end{bmatrix}$.
\end{itemize}
\end{example}

\begin{corollary}\label{order4}
Let $t_{k}$ be the number of elements of particular order $k$ in the Ap\'{e}ry set $Ap(\mathfrak{S}_{m},n_{1})$ are given by the following table.
\begin{table}[ht]
\caption{Number of elements of a particular order in the Ap\'{e}ry set}  
\centering 
\begin{tabular}{|c|c|c|c|c|c|c|c|c|} 
\hline\hline 
Order & $1$ & $2$ & $3$& $\cdots$ & $k$ & $\cdots$ & $m-1$& $m$\\ [0.5ex] 
\hline 
$A_{1}$ & $1$ & $1$ & $1$ &$\cdots$ & $1$ &$\cdots$ & $1$ & $1$\\\hline 
$A_{2}$ & $1$ & $1$ & $1$ &$\cdots$& $1$ & $\cdots$& $1$ & $0$\\\hline
$A_{3}$ & $1$ & $1$ & $1$ &$\cdots$& $1$ &$\cdots$& $1$ & $0$\\\hline
$A_{4}$ & $0$ & $1$ & $2$ &$\cdots$ & $k-1$ &$\cdots$ &$m-2$& $m-1$\\\hline
$A_{5}$ & $0$ & $1$ & $2$ &$\cdots$ &$k-1$ &$\cdots$ &$m-2$& $0$\\ \hline
$\textbf{Total}$ & $3$ & $5$ & $7$ &$\cdots$ & $2k+1$ & $\cdots$ & $2m-1$ & $m$\\ \hline 
\end{tabular} 
\end{table}
\end{corollary}
\proof The proof follows from the theorem \ref{aperytable4} and \ref{atars}.\qed

\begin{corollary}\label{exptangent} Let $I=\langle t^{n_{1}}\rangle$. The tangent cone $G_{\mathfrak{m}}(\mathfrak{S}_{m})$ of $\mathfrak{S}_{m}$ is a free $F(I)$-module. Moreover 
$$ G_{\mathfrak{m}}(\mathfrak{S}_{m})= \displaystyle\bigoplus_{k=1}^{m-1}(F(I)(-k))^{2m-1}\bigoplus(F(I)(-m))^{m}.$$
\end{corollary}
\proof The proof follows from the theorems \ref{tangentcone}, \ref{unique4}, \ref{aperytable4} and lemma \ref{order4}.\qed

\begin{corollary} The tangent cone $G_{\mathfrak{m}}(\mathfrak{S}_{m})$ is Cohen-Macaulay.
\end{corollary}
It is easily followed from the fact that $G_{\mathfrak{m}}(\mathfrak{S}_{m})$ is a free $F(I)$-module (see section 4 in \cite{cz2}). \qed

\begin{corollary} Let $HG_{\mathfrak{m}}(\mathfrak{S}_{m}(x))$ be the Hilbert series of $G_{\mathfrak{m}}$. Then $$HG_{\mathfrak{m}}(\mathfrak{S}_{m}(x))=\displaystyle \left(\sum_{k=1}^{m-1} (2k-1)x^{k}+mx^{m}\right)/(1-x).$$
\end{corollary}
\proof The proof follows from Corollary \ref{exptangent}.\qed
\medskip

\begin{remark}
Cohen-Macaulayness of the tangent cone of Arslan curves has been studied in \cite{ars}. Here we have calculated the Ap\'{e}ry table in detail to describe its tangent cone at the origin.
\end{remark}

\bibliographystyle{amsalpha}
                                                                          
\end{document}